
\magnification1200
\input amstex.tex
\documentstyle{amsppt}
\nopagenumbers
\hsize=12.5cm
\vsize=18cm
\hoffset=1cm
\voffset=2cm

\footline={\hss{\vbox to 2cm{\vfil\hbox{\rm\folio}}}\hss}

\def\DJ{\leavevmode\setbox0=\hbox{D}\kern0pt\rlap
{\kern.04em\raise.188\ht0\hbox{-}}D}

\baselineskip=13pt
\def\hf{{\textstyle{1\over2}}}
\def\a{\alpha}
\def\d{{\,\roman d}}
\def\e{\varepsilon}

\def\G{\Gamma}
\def\k{\kappa}

\def\t{\theta}
\def\={\;=\;}

\def\zt{\zeta(\hf+it)}

\def\D{\Delta}

\def\z{\zeta}

 \def\t{\theta}
\def\hf{{\textstyle{1\over2}}}

\def\le{\leqslant} \def\ge{\geqslant}
\font\tenmsb=msbm10
\font\sevenmsb=msbm7
\font\fivemsb=msbm5
\newfam\msbfam
\textfont\msbfam=\tenmsb
\scriptfont\msbfam=\sevenmsb
\scriptscriptfont\msbfam=\fivemsb
\def\Bbb#1{{\fam\msbfam #1}}

\def \NN {\Bbb N}

\font\ff=cmr8

\baselineskip=13pt

\font\teneufm=eufm10
\font\seveneufm=eufm7
\font\fiveeufm=eufm5
\newfam\eufmfam
\textfont\eufmfam=\teneufm
\scriptfont\eufmfam=\seveneufm
\scriptscriptfont\eufmfam=\fiveeufm
\def\mathfrak#1{{\fam\eufmfam\relax#1}}

\font\tenmsb=msbm10
\font\sevenmsb=msbm7
\font\fivemsb=msbm5
\newfam\msbfam
     \textfont\msbfam=\tenmsb
      \scriptfont\msbfam=\sevenmsb
      \scriptscriptfont\msbfam=\fivemsb
\def\Bbb#1{{\fam\msbfam #1}}

\def \NN {\Bbb N}

  \def\rightheadline{{\hfil{\ff
  On some upper bounds for $|\zt|$ and a divisor problem}\hfil\tenrm\folio}}

  \def\leftheadline{{\tenrm\folio\hfil{\ff
   Aleksandar Ivi\'c }\hfil}}
  \def\emptyheadline{\hfil}
  \headline{\ifnum\pageno=1 \emptyheadline\else
  \ifodd\pageno \rightheadline \else \leftheadline\fi\fi}

\topmatter
\title
ON SOME UPPER BOUNDS FOR THE ZETA-FUNCTION AND THE DIRICHLET DIVISOR  PROBLEM
\endtitle
\author   Aleksandar Ivi\'c
 \endauthor

\nopagenumbers

\medskip

\address
Aleksandar Ivi\'c, Katedra Matematike RGF-a
Universiteta u Beogradu, \DJ u\v sina 7, 11000 Beograd, Serbia
\endaddress

\bigskip
\email {
\tt
aleksandar.ivic\@rgf.bg.ac.rs, aivic\_2000\@yahoo.com }\endemail
\dedicatory
\enddedicatory
\abstract
{Let $d(n)$ be the number of divisors of $n$, let
$$
\D(x) := \sum_{n\le x}d(n) - x(\log x + 2\gamma -1)
$$
denote the error term in the classical Dirichlet
divisor problem, and let $\z(s)$ denote the Riemann
zeta-function. Several upper bounds for integrals of the type
$$
\int_0^T\D^k(t)|\zt|^{2m}\d t \qquad(k,m\in\Bbb N)
$$
are given. This complements the results of  the paper Ivi\'c-Zhai [9],
where asymptotic formulas for $2\le k \le 8,m =1$ were established
for the above integral.
 }
\endabstract
\endtopmatter

\document

\head
1. Introduction
\endhead

As usual, let
$$
\D(x) \;:=\; \sum_{n\le x}d(n) - x(\log x + 2\gamma - 1)\qquad(x\geqslant 2)
\leqno(1.1)
$$
denote the error term in the classical Dirichlet divisor problem (see e.g.,
Chapter 3 of [5]). Also let
$$
E(T) \;:=\;\int_0^T|\zt|^2\d t - T\Bigl(\log\bigl({T\over2\pi}\bigr) + 2\gamma - 1
\Bigr)\qquad(T\geqslant2)\leqno(1.2)
$$
denote the error term in the mean square formula for $|\zt|$.
Here $d(n)$ is the number of all positive divisors of
$n$, $\z(s)$ is the Riemann zeta-function, and $ \gamma = -\G'(1) = 0.577215\ldots\,$
is Euler's constant.  In [6] the author proved
several results involving the mean values of $\D(x), E(t)$ and
$$
\eqalign{
\D^*(x) :&= -\D(x)  + 2\D(2x) - \hf\D(4x)\cr&
= \hf\sum_{n\le4x}(-1)^nd(n) - x(\log x + 2\gamma - 1),\cr}
\leqno(1.3)
$$
which is the ``modified'' divisor function, introduced and studied
 by M. Jutila [10], [11].

In [9] the author and W. Zhai  studied the
moments $\int_0^T \D^k(t)|\zt|^2\d t$, where $k \in \Bbb N$ is fixed, to investigate
the interplay between the two fundamental functions $\D(t)$ and $|\zt|^2$.
It was proved that
$$
\int_0^T\D(t)|\zt|^2\d t \;\ll\; T(\log T)^{4},\leqno(1.4)
$$
and if $k$ is a fixed integer for which $2\le k \le 8$, then we have
$$
\int_1^{T}\Delta^{k}(t)|\zeta(\hf+it)|^2\d t=c_1(k)T^{1+\frac k4}\log T+
c_2(k)T^{1+\frac k4}+O_\e(T^{1+\frac k4-\eta_k+\varepsilon}),\leqno(1.5)
$$
 where
$c_1(k)$ and $c_2(k)$ are explicit constants, and where
$$
\eta_2=  \eta_3=  \eta_4=1/10,\; \eta_5=3/80,\; \eta_6=35/4742,\;
 \eta_7=17/6312,\; \eta_8=8/9433.
$$
It was also shown how the value of $\eta_2$ can be improved to $\eta_2 = 3/20$.
It may be well conjectured that the asymptotic formula (1.5) holds for integers
$k \ge 9$ as well (with some $\eta_k>0$), although this is beyond reach at present.
This is in tune with (1.2) and the classical conjecture that
$$
\int_0^T  \D^k(t)\d t =
C_kT^{1+k/4}
+ O_\e(T^{1+k/4-c(k)+\e})\leqno(1.6)
$$
holds with an explicit constant $C_k$ and some $c(k)>0$, when $k>1$ is a given
natural number. We note that (1.6) is at present known to hold for $2\le k \le 9$
(see W. Zhai [14]). In what concerns (1.4) it was conjectured in [9] that
one has
$$
\int_1^T\D(t)|\zt|^2\d t =
{T\over4}\Bigl(\log\frac{T}{2\pi} + 2\gamma - 1\Bigr)
+ O_\e(T^{3/4+\e}),
$$
however obtaining any asymptotic formula for the integral in (1.4)
is difficult. Here and later $\e$ denotes arbitrarily small positive constants,
not necessarily the same ones at each occurrence, while $f \ll_{a,b,\ldots} g$
(same as $f = O_{a,b,\cdots}(g)$) means that the implied constant depends on $a,b,\ldots$ .

\medskip
\head
2. Statement of results
\endhead
\medskip
A natural continuation of the previous investigations related to the integral in (1.5)
is the estimation of the more general integral
$$
\int_0^T\D^k(t)|\zt|^{2m}\d t \qquad(k,m\in\Bbb N),\leqno(2.1)
$$
where $k \ge1, m>1$. One would naturally want to obtain non-trivial upper bounds,
where by trivial we mean bounds coming from the use of the currently best known
upper bounds
$$
\Delta(x)\;\ll_\e\; x^{\theta+\varepsilon},\quad \theta=131/416=0.3149\ldots\leqno(2.2)
$$
and
$$
\zt \;\ll_\e\; |t|^{32/205+\e},\quad 32/205 = 0.15609\ldots \,.\leqno(2.3)
$$
We note that the exponent 32/205 has been recently improved to $53/342 = 0.15497\ldots\,$
in a forthcoming paper by J. Bourgain [2].
In the case of $m=2$ in (2.1) one would naturally wish to use results on the fourth moment
of $|\zt|$. We have (see Ivi\'c - Motohashi [7],[8] and Y. Motohashi [12])
$$
\int_0^T|\zt|^4 \d t = TQ_4(\log T) + E_2(T), \quad E_2(T) = O(T^{2/3}\log^8T),\leqno(2.4)
$$
where $Q_4(x)$ is an explicit polynomial of degree four in $x$ with
leading coefficient $1/(2\pi^2)$. We also have (here and later $C$ denotes positive generic
constants)
$$
\int_0^T E_2^2(t)\d t \;\ll\; T^2\log^CT \leqno(2.5)
$$
but neither (2.3) nor (2.4)-(2.5) are sufficiently strong to obtain non-trivial results regarding
(2.1) when $m=2$.

Our results are contained in the following
\medskip
THEOREM 1. {\it We have
$$
\eqalign{
\int_0^T \D(t)|\zt|^4\d t &\;\ll\; T^{41/32}\log^CT\qquad(41/32 = 1,28125),\cr
\int_0^T \D^2(t)|\zt|^4\d t &\;\ll\; T^{25/16}\log^CT\qquad(25/16 = 1,56250),\cr
\int_0^T \D^3(t)|\zt|^4\d t &\;\ll\; T^{59/32}\log^CT\qquad(59/32 = 1,84375),\cr
\int_0^T \D^4(t)|\zt|^4\d t &\;\ll\; T^{17/8}\log^CT\qquad(17/8 = 2,125).\cr}\leqno(2.6)
$$
}

\medskip
THEOREM 2. {\it We have}
$$
\eqalign{
\int_0^T \D(t)|\zt|^6\d t \;&\ll\; T^{49/32}\log^CT\qquad(49/32 = 1,53125),\cr
\int_0^T \D^2(t)|\zt|^6\d t \;&\ll\; T^{29/16}\log^CT\qquad(29/16 = 1,8125).\cr}
\leqno(2.7)
$$

\medskip
{\bf Remark 1}. The values of the constants $C$ in (2.6) and (2.7) are not
important, since the exponents are certainly not the best possible ones, as
will be discussed later. Indeed,
if one assumes the classical conjecture
$$
\D(x) \;\ll_\e\; x^{1/4+\e}\leqno(2.8)
$$
and the famous Lindel\"of Hypothesis
$$
\zt \;\ll_\e\;|t|^\e,\leqno(2.9)
$$
then trivially we have, for natural numbers $k > 1, m \ge 1$
$$
\int_0^T\D^k(t)|\zt|^{2m}\d t \;\ll_{\e,k,m}\; T^{1+k/4+\e} .\leqno(2.10)
$$
Proving (2.10) in full generality is not possible nowadays, since neither (2.8) nor (2.9) is
yet known to be true. It is classical that the Lindel\"of Hypothesis
follows from the Riemann Hypothesis (that all complex zeros of $\z(s)$ have real parts 1/2);
see e.g. [5, Chapter 1]. However, (2.8) does not seem to follow from any known hypotheses, and
the best known exponent $\theta=131/416=0.3149\ldots$ in (2.2) is very far from the conjectural
exponent $1/4+\e$.

The upper bound in (2.10) can be probably sharpened, at least for some values of $k$ and $m$,
to an asymptotic formula of the form
$$
\int_0^T\D^k(t)|\zt|^{2m}\d t = T^{1+k/4}Q_{m^2}(\log T) + O_{\e,k,m}(T^{1+k/4+\e-\rho_{k,m}})
$$
for some constant $\rho_{k,m} >0$, where $Q_{m^2}(x)$ is a polynomial of degree $m^2$, whose
coefficients depend on $k$ and $m$.
\medskip
{\bf Remark 2}.
The methods of proofs of the results allow one
to carry over the results of Theorem 1 and  Theorem 2  to the
integrals where $\D(t)$ is replaced by $\D(\a t)$ or $\D^*(\a t)$ for
any given $\a>0$. Here
$$
\eqalign{
\D^*(x) :&= -\D(x)  + 2\D(2x) - \hf\D(4x)\cr&
= \hf\sum_{n\le4x}(-1)^nd(n) - x(\log x + 2\gamma - 1),\cr}
$$
which is the ``modified''  function in the divisor problem.
In view of F.V. Atkinson's classical explicit formula (see [1] and Chapter 15 of [5]) for $E(T)$,
which shows analogies between $\D(x)$ and $E(T)$, it turns out that
$\D^*(x)$ is a better analogue of $E(T)$ than $\D(x)$ itself.

\medskip
{\bf Remark 3}. Finally, as in [9], we indicate two possible generalizations
of our results. Namely the results can be generalized if $\D(x)$ is replaced
either by $P(x) := \sum\limits_{n\le x }r(n) - \pi x$, or
$A^*(t) := \sum\limits_{n\le t}a(n)n^{{1-\k\over2}}$. As usual, $r(n) = \sum_{n=a^2+b^2}1$
denotes
the number of ways $n$ may be represented as a sum of two integer squares, and
$a(n)$ the $n$-th Fourier
coefficient of $\varphi(z)$, a
normalized eigenfunction of weight $\k$ for the Hecke operators $T(n)$, that is,
$a(1)=1  $ and $  T(n)\varphi=a(n)\varphi $ for every $n \in
\NN$.
\bigskip
\head
3. Proofs of the Theorems
\endhead
\bigskip
The ingredients in the proof are the asymptotic formula (1.5) (with $k=8$),
results on upper bounds for the moments of $|\zt|$, and H\"older's classical
inequality for integrals in the form
$$
\int_a^b f_1(x)\ldots f_r(x)\d x \le {\left(\int_a^b {f}^{p_1}_1(x)\d x\right)}^{1/p_1}\cdots
{\left(\int_a^b {f}^{p_r}_r(x)\d x\right)}^{1/p_r},\leqno(3.1)
$$
where $p_1, p_2, \ldots ,p_r >0$ and $f_1(x), f_2(x), \ldots, f_r(x)\ge 0$
are integrable functions in $[a,b]\;(a<b)$, and
$$
\frac{1}{p_1} + \frac{1}{p_2} + \ldots + \frac{1}{p_r} =1.
$$
\medskip
The case $r=2, p_1 = p_2 = 2$ is the standard Cauchy-Schwarz inequality
for integrals.
\medskip
For the moments of $|\zt|$ we use the bound (this is [5, Theorem 8.3])
$$
\int_0^T|\zt|^A\d t \;\ll\; T^{1 + (A-4)/8}\log^{C(A)}T\qquad(4 \le A \le 12).\leqno(3.2)
$$
\medskip
The value of the constant $C(A)$ in (3.2) can be given explicitly, but as mentioned, its value is not
important for our applications.

\medskip
To prove the first upper bound in (2.6) note that
$$
\eqalign{&
\int_0^T\D(t)|\zt|^4\d t = \int_0^T\D(t)|\zt|^{1/4}|\zt|^{15/4}\d t
\cr&
\ll \left(\int_0^T\D^8(t)|\zt|^2\d t\right)^{1/8}\left(\int_0^T|\zt|^{30/7}\d t\right)^{7/8}
\cr&
\ll (T^3\log T)^{1/8}\left(\int_0^T|\zt|^{4+\frac27}\d t\right)^{7/8}
\cr &\ll T^{41/32}\log^CT,\cr}\leqno(3.3)
$$
where $C$ denotes positive, generic constants as already mentioned. Here we used
H\"older's inequality (3.1), (1.5) with $k = 8$ and (3.2).

\medskip
{\bf Remark 4}.
It is readily checked that
the exponent (3.3) cannot be improved by using trivial estimation coming from the
bounds in (2.2) and (2.3).

\medskip
{\bf Remark 5}. The idea in proving (3.3), and other upper bounds of Theorem 1 as well,
is to use (1.5) with $k=8$. However, not the full asymptotic formula implied by (1.5)
is used, but just the  upper bound $T^3\log T$. There is a possibility to obtain small
improvements on all exponents in Theorem 1 and Theorem 2 as follows. First, recall (see
(2.2)) that there exists a constant $\t$ such that $1/4 \le \theta<1/3$ and
$$\Delta(x)\ll_\e x^{\theta+\varepsilon}, \ \ E(t)\ll_\e
t^{\theta+\varepsilon}.\leqno(3.4)
$$
In particular,  we can take $\theta=131/416=0.3149\cdots.$
The proofs of the bounds in (3.4) are due to M.N. Huxley [4] and N. Watt [13], respectively,
and they are the sharpest ones known.
Then for any $A$
satisfying $0\le A\le 11$ we have
$$\int_1^T |\Delta(x)|^A\d x\;\ll_\e\; T^{1+M(A)+\varepsilon}\leqno(3.5)
$$
and
$$
\int_1^T |E(t)|^A\d t  \;\ll_\e\; T^{1+M(A)+\varepsilon},\leqno(3.6)
$$
where
$$
M(A):=\max\left(\frac{A}{4}, \theta(A-2)\right).\leqno(3.7)
$$
This follows by the discussion given in the author's monograph [5, Chapter 13].

\medskip
For completeness, we also note that, for real $k \in [0,9]$, the limits
$$
E_k
\;:=\; \lim_{T\to\infty}T^{-1-k/4}\int_0^T|E(t)|^k\d t\leqno(3.8)
$$
exist. The analogous result holds also for the moments of $\D(t)$.
This was proved by
D.R. Heath-Brown [3], who used (3.5) and (3.6) in his proof. He also showed that
the limits of moments (both of $\D(t)$ and $E(t)$) without
absolute values also exist when  $k=1,3,5,7$ or 9. For the asymptotic
formulas for the moments of $\D(t), E(t)$ see W. Zhai [14], [15].
The merit of (3.8) that it gets rid of ``$\e$'' and establishes the existence
of the limit (but without an error term). Note that, with $\theta=131/416=0.3149\cdots,$
in (3.7) we have $M(A) = A/4$ for $A \le 262/27 = 9.\overline{703}$.
Using the method of [9] we can find a constant $8 < A_0 < A = 262/27$ for which
the bound
$$
\int_0^T|\D(t)|^{A_0}|\zt|^2\d t \;\ll_\e\; T^{1+A_0/4+\e}\leqno(3.9)
$$
will hold. Hence an improvement of
(3.3) will consist by using H\"older's inequality in such a way that instead of the
integral of $\D^8(t)|\zt|^2$ we have the integral of $|\D(t)|^{A_0}|\zt|^2$ with
$A_0$ as in (3.9). However, this would entail unwieldy exponents and the improvement would
not be large, so we worked out explicitly the results using only the integral of $\D^8(t)|\zt|^2$.

\medskip
We continue with the proof of the remaining bounds in (2.6). We have
$$
\eqalign{
&
\int_0^T\D^2(t)|\zt|^4\d t = \int_0^T\D^2(t)|\zt|^{1/2}|\zt|^{7/2}\d t
\cr&
\le \left(\int_0^T\D^8(t)|\zt|^2\d t\right)^{1/4}\left(\int_0^T|\zt|^{4+2/3}\d t\right)^{3/4}
\cr&
\ll (T^3\log T)^{1/4}(T^{1+1/12})^{3/4}\log^CT \ll T^{25/16}\log^CT,
\cr}
$$
$$
\eqalign{
&
\int_0^T\D^3(t)|\zt|^4\d t = \int_0^T\D^3(t)|\zt|^{3/4}|\zt|^{13/4}\d t
\cr&
\ll \left(\int_0^T\D^8(t)|\zt|^2\d t\right)^{3/8}\left(\int_0^T|\zt|^{4+6/5}\d t\right)^{5/8}
\cr&
\ll (T^3\log T)^{3/8}(T^{1+3/20})^{5/8}\log^CT \ll T^{59/32}\log^C T,
\cr}
$$
$$
\eqalign{
&
\int_0^T\D^4(t)|\zt|^4\d t = \int_0^T\D^2(t)|\zt||\zt|^{3}\d t
\cr&
\le \left(\int_0^T\D^8(t)|\zt|^2\d t\right)^{1/2}\left(\int_0^T|\zt|^{6}\d t\right)^{1/2}
\cr&
\ll (T^3\log T)^{1/2}(T^{1+1/4})^{1/2}\log^CT \ll T^{17/8}\log^CT.
\cr}
$$
This proves Theorem 1. The proof of Theorem 2 is on similar lines. Namely
$$
\eqalign{&
\int_0^T\D(t)|\zt|^6\d t = \int_0^T\D(t)|\zt|^{1/4}|\zt|^{23/4}\d t
\cr&
\ll \left(\int_0^T\D^8(t)|\zt|^2\d t\right)^{1/8}\left(\int_0^T|\zt|^{46/7}\d t\right)^{7/8}
\cr&
\ll (T^3\log T)^{1/8}\left(\int_0^T|\zt|^{4+\frac{18}{7}}\d t\right)^{7/8}
\cr &\ll T^{49/32}\log^CT,\cr}
$$
where again we used (1.5) with $k=8$, (3.1) and (3.2). Finally
$$
\eqalign{&
\int_0^T\D^2(t)|\zt|^6\d t = \int_0^T\D^2(t)|\zt|^{1/2}|\zt|^{11/2}\d t
\cr&
\le \left(\int_0^T\D^8(t)|\zt|^2\d t\right)^{1/4}\left(\int_0^T|\zt|^{22/3}\d t\right)^{3/4}
\cr &\ll (T^3\log T)^{1/4}(T^{1+10/24})^{3/4}\log^CT \ll T^{29/16}\log^CT,\cr}
$$
as asserted.

\medskip

\bigskip
\vfill
\eject
\topglue1cm
\medskip
\Refs
\medskip

\item{[1]} F.V. Atkinson, {\it The mean value of the Riemann zeta-function},
Acta Math. {\bf81}(1949), 353-376.

\item{[2]}  J. Bourgain, {\it Decoupling, exponential sums and the Riemann zeta-function},
preprint available at {\tt arXiv.1408.5794}.

\item{[3]} D.R. Heath-Brown,
{\it The distribution and moments of the error term in the Dirichlet
divisor problems}, Acta Arith. {\bf60}(1992), 389-415.

\item  {[4]} M. N. Huxley, {\it Exponential sums and lattice points III,}
Proc. London Math. Soc. {\bf87}(3)(2003), 591--609.

\item{[5]} A. Ivi\'c, {\it The Riemann zeta-function}, John Wiley \&
Sons, New York, 1985 (2nd ed. Dover, Mineola, New York, 2003).

\item{[6]} A. Ivi\'c, {\it On some mean value results for the zeta-function
and a divisor problem}, to appear in Filomat, preprint available
at {\tt arXiv:1406.0604}.

\item{[7]} A. Ivi\'c and Y. Motohashi, {\it
The mean square of the error term for the
    fourth moment of the zeta-function}, Proc. London Math. Soc. (3)
    {\bf69} (1994), 309-329.

\item{[8]} A. Ivi\'c and Y. Motohashi, {\it On the fourth power moment of the
Riemann zeta-function}, J. Number Theory {\bf51}(1995), 16-45.

\item{[9]} A. Ivi\'c and W.Zhai, {\it On some mean value results for $|\zt|$
and a divisor problem II}, Indagationes Mathematicae, 2015,
DOI 10.1016/j.indag.2015.05.002. Preprint available at {\tt arXiv:1502.00406}.

\item{[10]} M. Jutila, {\it Riemann's zeta-function and the divisor problem},
Arkiv Mat. {\bf21}(1983), 75-96 and II, ibid. {\bf31}(1993), 61-70.

\item{[11]} M. Jutila, {\it On a formula of Atkinson},
Topics in classical number theory, Colloq. Budapest 1981,
Vol. I, Colloq. Math. Soc. J\'anos Bolyai {\bf34}(1984), 807-823.

\item{[12]} Y. Motohashi, {\it Spectral theory of the Riemann zeta-function}, Cambridge
University Press, Cambridge, 1997.

\item{[13]} N. Watt, {\it A note on the mean square of $|\zt|$},
J. London Math. Soc. {\bf82(2)}(2010), 279-294.

\item{[14]} W. Zhai, {\it On higher-power moments of $\Delta(x)$}, Acta Arith.
{\bf112}(2004), 367-395; II. ibid. {\bf114}(2004), 35-54; III ibid. {\bf118}(2005),
263-281, and IV,
Acta Math. Sinica, Chin. Ser. {\bf49}\break (2006), 639-646.

\item{[15]} W. Zhai, {\it On higher-power moments of $E(t)$}, Acta Arith. {\bf115}(2004),
329-348.

\vfill

\endRefs

\enddocument

\bye